\documentclass[12pt]{amsart}
\usepackage{amsmath}
\usepackage{amssymb}
\usepackage{epsf}
\usepackage{epsfig}

\def\tp{\tilde p}

\def\F{{\bf F}}

\def\G{{\bf G}}
\def\H{{\bf H}}
\def\UU{{\bf U}}

\def\cal{\mathcal}

\def\M{{\bf M}}
\def\V{{\bf V}}

\def\C{{\mathbb C}}

\def\SG{{\mathfrak S}}

\def\std{{\rm std\,}}
\def\Des{{\rm Des\,}}

\def\tasse{{\rm C}}
\def\pigrec{{\underline{\pi}}}

\def\SS{{\bf S}}

\def\Sym{{\bf Sym}}     
\def\NCSF{{\bf Sym}}    
\def\FQSym{{\bf FQSym}} 
\def\WQSym{{\bf WQSym}} 
\def\FSym{{\bf FSym}}   
\def\QSym{{\it QSym}}   
\def\PBT{{\bf PBT}}     

\def\Pp{{\bf P}}
\def\Qq{{\bf Q}}

\def\LL{{\bf \Lambda}}
\def\WW{{\bf W}}

\newtheorem{example}{Example}[section]
\newtheorem{note}[example]{Note}

\newtheorem{theorem}[example]{Theorem}

\newtheorem{definition}[example]{Definition}
\newtheorem{proposition}[example]{Proposition}

\newtheorem{lemma}[example]{Lemma}
\def\mref#1{(\ref{#1})}

\def\up#1{\raise 1ex\hbox{\footnotesize#1}}

\def\Proof{\noindent \it Proof -- \rm}

\def\qed{\hspace{3.5mm} \hfill \vbox{\hrule height 3pt depth 2 pt width 2mm}
\bigskip}

\def\binomial#1#2{\left(\,\begin{matrix}#1 \cr #2\end{matrix}\,\right)}
\def\ncbinomial#1#2{\left[\,\begin{matrix}#1 \cr #2\end{matrix}\,\right]}

\def\U0{{\cal U}_0(gl_N)}

\def\carr{\blacktriangleright}

\def\<{\langle}
\def\>{\rangle}

\def\ra{\rightarrow}

\def\ie{{\it i.e.,\ }}

\def\vtr#1{\vrule height 0mm depth #1mm width 0mm}

\def\K{{\mathbb K}}

\def\shuff#1#2{\mathbin{
      \hbox{\vbox{
        \hbox{\vrule
              \hskip#2
              \vrule height#1 width 0pt
               }%
        \hrule}%
             \vbox{
        \hbox{\vrule
              \hskip#2
              \vrule height#1 width 0pt
               \vrule }%
        \hrule}%
}}}

\def\shuffl{{\mathchoice{\shuff{7pt}{3.5pt}}%
                        {\shuff{6pt}{3pt}}%
                        {\shuff{4pt}{2pt}}%
                        {\shuff{3pt}{1.5pt}}}}%
\def\shuffle{\, \shuffl \,}

\def\DessinMatrix#1{\vcenter{\hbox{\makebox[1.7ex]{$#1$}}}}
\def\GenMatrix#1{\vcenter{\halign{&$\DessinMatrix{##}$\cr#1}}\egroup}

\def\setinterlineskip#1{\baselineskip=0pt
  \lineskip=#1 \lineskiplimit=\maxdimen}

\def\matrice{%
  \bgroup
  \let\ =\omit
  \let\\=\cr
  \setinterlineskip{4.0pt}
  \GenMatrix}

\def\DessinsMatrix#1{\vcenter{\hbox{\makebox[1.3ex]{$\scriptstyle#1$}}}}
\def\GensMatrix#1{\vcenter{\halign{&$\DessinsMatrix{##}$\cr#1}}\egroup}

\def\smallmatrice{%
  \bgroup
  \let\ =\omit
  \let\\=\cr
  \setinterlineskip{3.0pt}
  \GensMatrix}

\newlength{\Hackl}
\newcommand{\Hack}{\vrule height \Hackl width 0pt}
\newcommand{\indexmat}%
    {\smallmatrice{\Hack a\\\Hack b\\\Hack c\\\Hack d\\\Hack e\\}}

\newdimen\Squaresize \Squaresize=14pt
\newdimen\Thickness \Thickness=0.5pt

\def\Square#1{\hbox{\vrule width \Thickness
   \vbox to \Squaresize{\hrule height \Thickness\vss
      \hbox to \Squaresize{\hss#1\hss}
   \vss\hrule height\Thickness}
\unskip\vrule width \Thickness}
\kern-\Thickness}

\def\Vsquare#1{\vbox{\Square{$#1$}}\kern-\Thickness}

\def\id{{\rm Id}}

\def\Sylv{{\rm Sylv}}   
\def\A{{\mathcal A}}
\def\cat{{\rm cat}}
\def\SS{{\bf S}}

\title[Noncommutative Symmetric Functions VII]%
{Noncommutative Symmetric Functions VII: \\
Free Quasi-Symmetric Functions Revisited}

\author[G.~H.~E. Duchamp, F. Hivert, J.-C. Novelli, J.-Y. Thibon]{%
G\'erard H.~E. Duchamp, Florent Hivert, Jean-Christophe Novelli, %
Jean-Yves Thibon}

\address[Duchamp]{Institut Galil\'ee, LIPN, CNRS UMR 7030\\
99, avenue J.-B. Clement, F-93430 Villetaneuse, France}
\address[Hivert]{LITIS, Universit\'e de Rouen ; Avenue de l'universit\'e ;
76801 Saint \'Etienne du Rouvray, France\\}
\address[Novelli, Thibon]{Institut Gaspard-Monge, Universit\'e Paris-Est, 
5, boulevard Descartes \\Champs-sur-Marne \\77454 Mar\-ne-la-Vall\'ee cedex 2 \\
France}
\email[G. H. E. Duchamp]{ghed@lipn.univ-paris13.fr}
\email[F. Hivert]{hivert@univ-rouen.fr}
\email[J.-C. Novelli]{novelli@univ-mlv.fr (corresponding author)}
\email[J.-Y. Thibon]{jyt@univ-mlv.fr}

\begin{document}

\begin{abstract}
We prove a Cauchy identity for free quasi-symmetric functions and apply
it to the study of various bases.
A free Weyl formula and a generalization of the splitting formula are also
discussed.
\end{abstract}

\maketitle

\section{Introduction}

This article is essentially an appendix to~\cite{NCSF6}.
We gather here some useful properties of the algebra $\FQSym$ of free
quasi-symmetric functions which were overlooked in~\cite{NCSF6}. 

Recall that $\FQSym$ is a subalgebra of the algebra of noncommutative
polynomials in infinitely many variables $a_i$ which is mapped onto Gessel's
algebra of quasi-symmetric functions $QSym$ by the commutative image
$a_i\mapsto x_i$ of $\K\<A\>$.
As an abstract algebra, it is isomorphic to the convolution algebra of
permutations introduced by Reutenauer and his school \cite{Re,MR,PR}, and
further studied in \cite{NCSF6,LR1,LR2,AS}.
However, the realization in terms of the variables $a_i$ provides one with a
better understanding of several aspects of the theory.
For example, it becomes possible, and sometimes straigthforward, to imitate
various constructions of the theory of symmetric (or quasi-symmetric)
functions. An illustration is provided by the construction of the
coproduct given in \cite{NCSF6}: a free quasi-symmetric function $F$
can be regarded as a ``function'' of a linearly ordered alphabet $A$, and the
obvious analog of the coproduct of $QSym$, that is,
$F(A)\mapsto F(A'\oplus A'')$, where $A'$ and $A''$ are two mutually commuting
copies of $A$, and $\oplus$ is the ordered sum, gives back the
coproduct of~\cite{MR}.

In the present text, we  further investigate the r\^ole of the auxiliary
variables $a_i$.

We start with an alternative definition of the standard basis of $\FQSym$,
as resulting from a noncommutative lift of a Weyl-type character formula.

Next, we formulate a \emph{free Cauchy identity} in $\FQSym$, and
investigate its implications. 
In the classical theory of symmetric functions, it is mostly in the so-called
{\it Cauchy identities}, involving the product $XY$ of two alphabets $X,Y$,
that the existence of the underlying variables manifests itself.
Even so, it is only when the auxiliary alphabet is \emph{specialized} that one
really sees the variables. Otherwise, the transformation $f(X)\mapsto f(XY)$
is best interpreted as a coproduct, dual to the internal product.
The present version involves two alphabets $A'$, $A''$, and specializes to the
noncommutative Cauchy formula of~\cite{GKLLRT} under $A'\mapsto X$,
$A''\mapsto A$.
This allows us to compute the commutative specializations of various bases
of $\FQSym$, in particular the basis $\V_\sigma$ of~\cite{NCSF6} and the basis
$\M_\sigma$ of~\cite{AS}, recovering here the result of~\cite{AS} in a simpler
way.

We first illustrate the free Cauchy identity by computing the Hopf duals of
various subalgebras of $\FQSym$ is a unified and straightforward way.
We then apply it to the study of several multiplicative bases of $\FQSym$ and
of their dual bases, including the free monomial functions of Aguiar and
Sottile~\cite{AS}, and the basis $\V_\sigma$ introduced in~\cite{NCSF6}.
We conclude with an extension of the splitting formula to $\FQSym$.

{\footnotesize
{\it Acknowledgements.}
This project has been partially supported by the grant ANR-06-BLAN-0380.
The authors would also like to thank the contributors of the MuPAD project,
and especially of the combinat part, for providing the development environment
for their research (see~\cite{HTm} for an introduction to MuPAD-Combinat).
}

\section{Background and notations}

Our notations are as in~\cite{NCSF6}. 
All algebras are over some field $\K$ of characteristic $0$. We shall also
need the following notations.

\subsection{}
Recall that for two permutations
$\sigma$ and $\tau$, $\sigma\bullet\tau$ denotes their shifted concatenation. 
In $\FQSym$, we set
\begin{equation}
\G_\sigma\bullet \G_\tau=\G_{\sigma\bullet\tau}
\end{equation}
so that
$I=(i_1,\dots,i_r)$ and if $\alpha_k\in\SG_{i_k}$ for all $k\in[1,r]$, then
\begin{equation}
\label{prodGsI}
 \G_{\alpha_1} \dots \G_{\alpha_r} =
(\G_{\alpha_1}\bullet \dots\bullet \G_{\alpha_r}) * S^I.
\end{equation}
Indeed, remembering that the internal product in the $\G$-basis is
opposite to the product in the symmetric group ($\G_\sigma*
\G_\tau=\G_{\tau\sigma}$), we recognize the convolution product of
permutations.

\subsection{}
It has been shown in \cite{FH} that quasi-symmetric polynomials could be
defined as the invariants of a certain action of the symmetric group, or of
the Hecke algebra, called the quasi-symmetrizing action. This action
will be denoted by underlining the elements, \emph{e.g.},
$\underline{\sigma}\cdot f$ denotes the quasi-symmetrizing action of a
permutation $\sigma$ on a polynomial $f$.

\section{The free Weyl formula}

As already mentioned, $\QSym(X)$ is the algebra of $\SG[X]$-invariants in
$\K[X]$ for a special action of the symmetric group, called the
quasi-symmetrizing action.
Moreover, the analogy between Gessel's functions $F_I$ and Schur functions has
been traced back to the existence of an expression of $F_I$ by a kind of Weyl
character formula involving this action.
In this Section, we show that this Weyl formula can be lifted to $\FQSym$.

For a word $w$ on $\{a_1,\dots,a_n\}$, let $|w|_i$ be the number of
occurrences of the letter $a_i$ in $w$ and let the monomial
$m=x^{|w|}=x_1^{|w|_1}x_2^{|w|_2} \dots x_n^{|w|_n}$ be the commutative image
of $w$. We denote by $\tasse(m)$ the composition obtained by removing the
zeros in the vector $|w|$.
The construction relies on the following well-known lemma
(see~\cite{NCSF4,Nov}):
\begin{lemma}
The map $\phi:w \mapsto (x^{|w|}, \std(w)^{-1})$ is a bijection between words
and pairs $(m,\sigma)$ such that $\tasse(m)$ is finer than $C(\sigma)$.
\end{lemma}

Hence, for each $\sigma$, there is a linear isomorphism $\phi_\sigma$
between the subspace of $\K\<\<A\>\>$ spanned by the words of standardized
$\sigma$ and the subspace of $\K[X]$ spanned by the monomials such that
$I(m)\leq C(\sigma)$.

The definitions of $F_I$ and $\F_\sigma$ can then be rewritten as
\begin{equation}
\label{defFrewrite}
F_I = \sum_{I(m) \leq I} m\,, \qquad\text{and}\qquad
\F_\sigma = \sum_{I(m) \leq C(\sigma)} \phi^{-1} (m, \sigma)
= \phi_\sigma^{-1}(F_I)\,.
\end{equation}

Recall that the degenerate Hecke algebra $H_n(0)$ is generated by elements
$(\pi_1,\ldots,\pi_{n-1})$ satisfying the braid relation and $\pi_i^2=\pi_i$.
The Schur functions can be obtained in terms of an action of $H_n(0)$ on
$\C[X]$ defined by
\begin{equation}
\pi_i(f) := \frac{x_if - \sigma_i(x_if)}{x_i-x_{i+1}}
\end{equation}
where $\sigma_i$ is the automorphism exchanging $x_i$ and $x_{i+1}$. In this
case, $s_\lambda=\pi_{\omega}(x^\lambda)$, where $\omega=(n,(n-1),\ldots,1)$
is the longest element of $\SG_n$. If one instead uses the quasi-symmetrizing
action, denoted here as in~\cite{FH} by
$f\mapsto \underline{\sigma}_i f$,
one still has an action of $H_n(0)$ and by complete symmetrization of
monomials, one obtains Gessel's functions
\begin{equation}
F_I = \pigrec_{\omega}(x^I)\,.
\end{equation}
From the definition of~\cite{FH}, the monomials $m'$ appearing in
$\pigrec_\omega(x^I)$ satisfy $\tasse(m')<I$. Hence
it makes sense to define an action of $H_n(0)$ on $\K\<A\>$ by lifting the
action on $\K[X]$ along the maps $\phi_\sigma$.

\begin{definition}
Let $f\in H_n(0)$ and $w\in A^*$. The action of $f$ on $w$ is defined by
\begin{equation}
\underline{f}(w) = \phi_{\std(w)^{-1}}^{-1} (\underline{f}(|w|)).
\end{equation}
where $\underline{f}(|w|)$ is the quasi-symmetrizing action.
\end{definition}

This action can be computed as follows. For a word $w$ and $k<|w|_i$, denote
by $w(a_i\uparrow k)$ (resp. $w(a_i\downarrow k)$), the word obtained by
replacing the last (resp. first) $k$ letters $a_i$ by $a_{i+1}$ (resp.
$a_{i-1}$).
Remark that
\begin{equation}
\std(w(a_i\uparrow k)) = \std(w(a_i\downarrow k)) = \std(w)\,,
\end{equation}
so that the action of $\pigrec_i$ on $w$ is given by
\begin{equation}
\pigrec_i(w) =\left\{
\begin{array}{cl}
  w & \text{if $|w|_i = 0$ and $|w|_{i+1} = 0$,}\\
  w & \text{if $|w|_i \neq 0$ and $|w|_{i+1} \neq 0$,}\\
  \sum\limits_{k=0}^{|w|_i} w(a_i\uparrow k)
    & \text{if $|w|_i \neq 0$ and $|w|_{i+1} = 0$,}\\
  - \sum\limits_{k=1}^{|w|_i-1} w(a_i\downarrow k)
    & \text{if $|w|_i = 0$ and $|w|_{i+1} \neq 0$.}\\
\end{array}\right.
\end{equation}

Then the quasi-symmetric Weyl formula together with
Equation~(\ref{defFrewrite}) gives:
\begin{theorem}
Let $\sigma\in\SG_n$ and $A=\{a_1,\dots a_n\}$. Then
\begin{equation}
\F_\sigma(A) = \pigrec_\omega (A_\sigma)\,,
\end{equation}
where $A_\sigma$ is the smallest word on $A$ of standardized word
$\sigma^{-1}$ for the lexicographic order.
\end{theorem}

Note that $A_\sigma$ is the unique word with standardized $\sigma^{-1}$ and
evaluation $C(\sigma)$. For example, if $\sigma=1472635$ then
$\sigma^{-1}=1462753$ and $A_\sigma=a_1a_2a_3a_1a_3a_2a_1$. 

The full algebra of invariants of the free quasi-symmetrizing action
(either of $\SG_n$ or of $H_n(0)$) is denoted by $\WQSym(A)$.
In the limit $n\rightarrow\infty$, it acquires the structure of a self-dual
Hopf algebra~\cite{FH,NT1,HNT08}.
Its bases are parametrized by set compositions (ordered set partitions), and
the dimensions of its graded components are given by the ordered Bell numbers.

\section{The free Cauchy identity}

One of the most basic formulas in the classical theory of symmetric functions
is the so-called \emph{Cauchy identity}, which can be interpreted as giving an
expression of the identity map of $Sym$, identified to an element of
$Sym\otimes Sym \simeq Sym(X,Y)$, in terms of the variables $x_i,y_j$.

In the commutative case, the Cauchy identity encodes almost entirely the Hopf
algebra structure. The Cauchy kernel
\begin{equation}
K(X,Y)=\sigma_1(XY)=\prod_{i,j}(1-x_iy_j)^{-1}
\end{equation}
is reproducing for the standard scalar product $\<\,,\,\>$ defined by
\begin{equation}
\<s_\lambda(X),s_\mu(X)\>_X=\delta_{\lambda\mu}\,,
\end{equation}
that is,
\begin{equation}
f(X)=\<K(X,Y), f(Y)\>_Y\,.
\end{equation}
Also, the antipode of $Sym$ is the map $\tilde\omega: f(X)\mapsto f(-X)$ is
given by
\begin{equation}
\tilde\omega f(X)=f(-X)=\<K(-X,Y),f(Y)\>_Y= \<\lambda_{-1}(XY),f(Y)\>_Y\,,
\end{equation}
that is, its kernel is the inverse (in the sense of formal series in $X$, $Y$)
of the Cauchy kernel. Similarly, the Adams operations in the sense
of Hopf algebras (the convolution powers of the identity map) have
as kernels the corresponding powers $\sigma_1(kXY)$ of the Cauchy
kernel. Finally, the coproduct itself is given by
\begin{equation}
\Delta (f)= f(X+Y) = D_{K(X,Y)} f(X)\,,
\end{equation}
where $ D_{K(X,Y)}$ is the adjoint of the operator of the multiplication
by $K(X,Y)$, with respect to $\<\,,\>_X$.

\medskip
The interpretation of $\FQSym$ as a subalgebra of $\K\<A\>$ allows one to lift
the Cauchy identity to free quasi-symmetric functions.

Let $x_{ij}=$
\scriptsize$\binomial{i}{j}$\normalsize\
be commuting indeterminates, and $a_{ij}=$
\scriptsize$\ncbinomial{i}{j}$\normalsize\
be noncommuting ones. We shall denote by
\scriptsize $\binomial{i_1\ i_2 \cdots i_r}{j_1\ j_2 \cdots j_r}$ \normalsize\
the monomial
\scriptsize $\binomial{i_1}{j_1}\binomial{i_2}{j_2}\cdots \binomial{i_r}{j_r}$ 
\normalsize\
and by
\scriptsize $\ncbinomial{i_1,i_2,\cdots i_r}{j_1,j_2,\cdots j_r}$ \normalsize\
the word
\scriptsize\
$\ncbinomial{i_1}{j_1}\ncbinomial{i_2}{j_2}\cdots \ncbinomial{i_r}{j_r}$. 
\normalsize\
Such expressions will be referred to respectively as \emph{bimonomials} and
\emph{biwords}. With a bimonomial
\scriptsize $\binomial{u}{v}$\normalsize,\
we associate the two biwords
\scriptsize$\ncbinomial{u'}{v'}$\normalsize\
and
\scriptsize$\ncbinomial{u''}{v''}$\normalsize\
obtained by sorting
\scriptsize$\binomial{u}{v}$\normalsize\
in lexicographic order with priority respectively on the top and on the bottom
row. For example, with
\begin{equation}
\binomial{u}{v}=\binomial{cabaabcba}{abccaacbb}\,,
\end{equation}
we get
\begin{equation}
\ncbinomial{u'}{v'}=\binomial{aaaabbbcc}{abbcabcac} \quad\text{and}\quad
\ncbinomial{u''}{v''}=\ncbinomial{abcaababc}{aaabbbccc}.
\end{equation}
Observe that $\std(v')=145726839=(158236479)^{-1}=\std(u'')^{-1}$, which, in
fact, is a general property.

\begin{lemma}
With the previous notations, for all bimonomials
\scriptsize$\binomial{u}{v}$\normalsize, one has $\std(u'')=\std(v')^{-1}$ and
all pairs of words with this property are obtained by this process. Therefore,
we have an evaluation-preserving bijection between such pairs of words and
bimonomials.
\end{lemma}

\Proof
If $u$ and $v$ are permutations, cleary
$u''\otimes v' = \sigma^{-1}\otimes\sigma$ where $\sigma$ is the permutation
such that
\begin{equation}
\binomial{u}{v} = \binomial{\id}{\sigma} = \binomial{\sigma^{-1}}{\id}.
\end{equation}
Now, the bi-standardization process satisfies:
\begin{equation}
\left\<\binomial{\std u}{\std v}\right>
= \std u'' \otimes \std v',
\end{equation}
whence the result.
\qed

Let $\<\,\>$ denote the linear map from $\K[[x_{ij}]]$ to
$\K\<\<A\>\>\otimes \K\<\<A\>\>$ defined by
$\<\,$%
\scriptsize$\binomial{u}{v}$\normalsize
$\,\>=u''\otimes v'$. From the lemma, we obtain:

\begin{proposition}
For any pair of adjoint bases $({\bf P},{\bf Q})$ of $\FQSym$, one has
\begin{equation}
\left\<\prod_{i,j}\frac{1}{1-x_{ij}}\right\>=
\sum_{\std(v)=\std(u)^{-1}}u\otimes v=
\sum_{\sigma} \F_{\sigma}\otimes \G_{\sigma} =
\sum_{\sigma} {\bf P}_{\sigma}\otimes {\bf Q}_{\sigma}.
\end{equation}
\end{proposition}

In particular, identifying the tensor product $\FQSym\otimes \FQSym$ with
$\FQSym(A',A'')$, where $A'$ and $A''$ are two commuting alphabets, and
specializing $a'_i\mapsto x_i$ (commuting variables), $a''_i\mapsto a_i$, we
get
\begin{equation}
\sum_{\sigma}\F_{\sigma}(X)\G_{\sigma}(A)=\sum_{I}F_I(X)R_I(A)=
\prod_{i,j}^{\rightarrow}(1-x_ia_j)^{-1}=\sigma_1(XA)\,,
\end{equation}
and for any pair $({\bf P},{\bf Q})$ of adjoint bases, the evaluation morphism
\begin{equation}
ev_{X,A} :u\otimes v\mapsto u(X) v(A)
\end{equation}
maps $\sum_{\sigma}{\bf P}_{\sigma}\otimes {\bf Q}_\sigma$ to $\sigma_1(XA)$.

\section{Application to some duality questions}

\subsection{Specialization of the free Cauchy formula to $\NCSF$}

From the free Cauchy formula, one derives a new way to recover the Cauchy
identity associated with the dual pair of Hopf algebras $(\QSym,\Sym)$.
Instead of mapping the $a_i'$ to our commuting variables $x_i$, let us set
\begin{equation}
a_i'\mapsto \tilde a_i \,,
\end{equation}
where the $\tilde a_i$'s satisfy the hypoplactic relations. Let $\equiv_H$
denote the hypoplactic congruence on words on $A$ as well as on permutations
considered as words on the positive integers.
Recall that $u\equiv_H v$ iff (1) they have the same evaluation, and
(2) $\sigma=\std(u)^{-1}$ and $\tau=\std(v)^{-1}$ have the same descent set
(see~\cite{NCSF4,Nov}). This implies that
\begin{equation}
\sigma\equiv_H\tau \ \Longleftrightarrow\  \G_\sigma (A) \equiv_H \G_\tau(A)\,,
\end{equation}
and
\begin{equation}
\F_\sigma(A) \equiv_H \F_\tau(A) \ \Longleftrightarrow C(\sigma)=C(\tau)\,.
\end{equation}
Denote by $\tilde w$ the image of $w$ in the hypoplactic algebra.
If $C(\sigma)=I$, then $\tilde \F_\sigma=F_I(\tilde A)$, the hypoplactic
quasi-symmetric function,
so that
\begin{equation}
\begin{split}
\sum_{\std(u)=\std(v)^{-1}}\tilde u\otimes v
& =\sum_\sigma \tilde \F_\sigma\otimes \G_\sigma
  =\sum_I F_I(\tilde A)\otimes \left(\sum_{C(\sigma)=I}\G_\sigma\right) \\
& =\sum_IF_I(\tilde A)\otimes R_I(A)\,.
\end{split}
\end{equation}
This formula reflects the fact that the dual $\Sym^*=\QSym$ of $\Sym$ can
be identified with the quotient of $\FQSym$ by the hypoplactic congruence.

\subsection{Specialization of the free Cauchy formula to $\FSym$}

This is reminescent of a result of Poirier-Reutenauer~\cite{PR} allowing to
identify the dual of $\FSym$ to the quotient of $\FQSym$ by the relations
\begin{equation}
\G_\sigma \sim \G_\tau\ \Longleftrightarrow \ \sigma \equiv_P \tau\,,
\end{equation}
where $\equiv_P$ is the plactic congruence, permutations being considered
as words. In fact, two words $u,v\in A^*$ are plactically equivalent
iff, as above, (1) they have the same evaluation, and (2)
$\std(u) \equiv_P \std(v)$. Hence,
\begin{equation}
\G_\sigma \sim \G_\tau \ \Longleftrightarrow \
\G_\sigma(A) \equiv_P \G_\tau(A)\,,
\end{equation}
the plactic congruence being understood in the free associative algebra on
$A$. This shows that $\FSym^*$ is actually a subalgebra of the plactic
algebra, spanned by the plactic classes $\bar\G_\sigma$ of the $\G_\sigma$'s.
If $P(\sigma)=t$, let us set $\bar\G_\sigma=\SS_t^*$ (we denote by $\bar w$
the image of $w$ in the plactic algebra), and consider the specialization
\begin{equation}
\begin{split}
\sum_{\std(u)=\std(v)^{-1}}\ u\otimes \bar v
& =\sum_\sigma  \F_\sigma\otimes \bar \G_\sigma  
  =\sum_t\left(\sum_{P(\sigma)=t}\F_\sigma\right)\otimes\SS_t^*\\
& =\sum_t\SS_t\otimes \SS_t^*\,.
\end{split}
\end{equation}

This last equality looks closer to the classical Cauchy identity. However, it
is not symmetric due to the non-self-duality of $\FSym$. To recover a
symmetric form, let us now impose the plactic relations on the first factor of
the tensor product (so that the resulting element belongs to
${\rm Pl(A)}\otimes {\rm Pl}(A)$). We get
\begin{equation}
\sum_t \bar \SS_t\otimes \SS_t^*=\sum_\lambda S_\lambda\otimes S_\lambda\,,
\end{equation}
where $S_\lambda$ is the plactic Schur function
\begin{equation}
S_\lambda=\sum_{{\rm shape}(\bar w)=\lambda}\bar w\,,
\end{equation}
defined as the sum of all tableaux of shape $\lambda$ in the plactic algebra.
We see that
\begin{equation}
S_\lambda =\sum_{{\rm shape}(t)=\lambda}\SS_t^*\,,
\end{equation}
that is, the commutative algebra spanned by the plactic Schur functions
introduced by Lascoux and Sch\"utzenberger~\cite{LS} gets naturally identified
with a subalgebra of $\FSym^*$.

\subsection{Specialization of the free Cauchy formula to $\PBT$}

By the same reasoning, we see that the dual $\PBT^*$ of the Loday-Ronco
algebra of planar binary trees~\cite{LR1,HNT05} $\PBT$ is ismomorphic to the
image of $\FQSym$ under the canonical projection
\begin{equation}
\pi :\ \K\<A\>\longrightarrow \K[\Sylv(A)]\simeq \K\<A\>/\equiv\,,
\end{equation}
where $\Sylv(A)$ denotes the sylvester monoid.
The dual basis $\Qq_T$ of the natural basis $\Pp_T$ is $\Qq_T=\pi(\G_{w_T})$,
where $w_T$ is the permutation canonically associated with the tree $T$.

\section{Multiplicative bases and their adjoint bases}

The free Cauchy formula is also useful for studying duals of multiplicative
bases, from which bases of primitive elements can be extracted. In this
section, we discuss two kinds of analogs of the monomial symmetric functions
and of the forgotten symmetric functions. One of them is the free monomial
basis of Aguiar-Sottile~\cite{AS}.

\subsection{Analog of the pair $(h,m)$ (first kind)}

We know that $\Sym\subset \FQSym$, and that the noncommutative products of
complete homogenous functions $S^I$ are realized by setting
$S_m=\G_{(12\ldots m)}$ and $S^I=S_{i_1}S_{i_2}\cdots S_{i_r}$, as usual.
There are of course many possibilities to extend $S^I$ to a multiplicative
basis of $\FQSym$.
We know that $\FQSym$ is free as an algebra, and that its
number of algebraic generators by degree is given by Sequence~A003319
of~\cite{Sloane}, counting \emph{connected permutations}. In the $\G$ basis,
we make use of \emph{anticonnected} permutations, \ie permutations $\sigma$
whose mirror image $\bar\sigma$ is connected.
If we introduce the \emph{left-shifted concatenation} of words
\begin{equation}
u\carr v = u[l]v \quad {\rm if}\ u\in A^k,\ v\in A^l\,,
\end{equation}
similar to $\bullet$, but with the shift on the first factor, we can start
with the convention
\begin{equation}
\label{eq:multS}
\SS^\sigma =\SS^{\sigma_1}\SS^{\sigma_2}\cdots \SS^{\sigma_r}
\end{equation}
whenever $\sigma=\sigma_1\carr\sigma_2\carr\cdots\carr\sigma_r$ is the
decomposition of $\sigma$ into anticonnected components. It remains to decide
the value of $\SS^\sigma$ when $\sigma$ is anticonnected. For $\sigma=id_m$,
we want $\SS^\sigma=S_m=\G_\sigma$, but for other anticonnected permutations,
there are several possibilities. Note that, with our indexation scheme, we
have $S^I=\SS^{\omega(I)}$, so that
\begin{equation}
\SS^{\omega(I)}=\sum_{\tau\le\omega(I)}\G_\tau\,,
\end{equation}
where $\leq$ is the left weak order.
One possibility is to keep the rule $\SS^\sigma=\G_\sigma$, but the previous
equation suggests that a reasonable choice would be
\begin{equation}
\label{eq:defS}
\SS^\sigma := \sum_{\tau\le\sigma}\G_\tau\,
\end{equation}
for $\sigma$ anticonnected. Then the resulting basis has the following
properties.

\begin{proposition}
(i) For any permutation $\sigma$ (anticonnected or not), $\SS^\sigma$ is given 
by formula \mref{eq:defS}.

\noindent
(ii) The adjoint basis of $(\SS^\sigma)$ is the free monomial
basis $(\M_\sigma)$ of Aguiar and Sottile \cite{AS}.
\end{proposition}

\Proof
(i) One has only to prove that $\SS^{\sigma\carr\tau}$ satisfies
formula~(\ref{eq:defS}), under the assumptions that both $\SS^{\sigma}$ and
$\SS^{\tau}$ satisfy it for all permutations $\sigma$ and $\tau$.
Let $n=|\sigma|$ and $m=|\tau|$.
It is first obvious that all elements in the expansion of
$\SS^{\sigma\carr\tau}$ are smaller than $\sigma\carr\tau$.

Conversely, let us consider a permutation $\mu$ smaller than $\sigma\carr\tau$.
Then let $\sigma'$ be the permutation obtained by standardizing the subword
consisting of letters of $\mu$ on the interval $[m+1,m+n]$. By definition of
the weak left order, $\sigma'\le\sigma$. Defining in the same way $\tau'$, and
checking that $\tau'\le\tau$, one then concludes that
$\mu\in (\sigma'[m]\shuffle\tau')$, so that $G_\mu$ appears in the expansion
$\SS^{\sigma\carr\tau}$.

(ii) Straightforward from Formula~1.12 p.~232 in~\cite{AS}.
\qed

For example, with $I=(1,2,1)$, we have
\begin{equation}
\begin{split}
\SS^{4132}  &= \G_{1} (\G_{123}+\G_{132}) \\
&= \G_{1234} + \G_{2134} + \G_{3124} + \G_{4123} \\
 &+ \G_{1243} + \G_{2143} + \G_{3142} + \G_{4132}.
\end{split}
\end{equation}
And indeed, the permutations smaller than, or equal to $4132$ in the left
permutohedron are the eight permutations listed above.
We also have
\begin{equation}
\begin{split}
\SS^{4231} =  S^{121}  &= \G_{1} \G_{12} \G_{1} \\
&= \sum_{\sigma\in 1\shuffle 23 \shuffle 4} \G_{\sigma} \\
&= \G_{1234} + \G_{2134} + \G_{3124} + \G_{4123}
 + \G_{1243} + \G_{2143}\\  &+ \G_{3142} + \G_{4132}
 + \G_{1342} + \G_{2341} + \G_{3241} + \G_{4231}.
\end{split}
\end{equation}
And indeed, the permutations smaller than, or equal to $4231$ in the left
permutohedron are the twelve permutations listed above.

As a consequence, we can easily derive the specialization $\M_\sigma(X)$ of
the free monomial basis to commuting variables, as given in \cite{AS}.

Let us recall that permutations $\tau$ whose descent composition $C(\tau)$
is $I$ form an interval $[\alpha(I),\omega(I)]$ of the weak order (left
permutohedron). We shall say that permutations that are $\omega(I)$ for an $I$
are \emph{bottom} descent permutations.

\begin{lemma}[Theorem~7.3 of~\cite{AS}]
Let $\sigma$ be a permutation.
Then
\begin{equation}
\M_\sigma(X) =
\left\{
\begin{array}{cc}
M_I &\text{if $\sigma=\omega(I)$}, \\
0   &\text{otherwise}, \\
\end{array}
\right.
\end{equation}
\end{lemma}

\Proof
Indeed, the free Cauchy formula tells us that on the one hand
\begin{equation}
\sum_\sigma \M_\sigma(X) \SS^\sigma(A) = \sum_I F_I(X) R_I(A)
=\sum_I F_I(X)\sum_{C(\tau)=I}\G_\tau(A)\,,
\end{equation}
and on the other hand,
\begin{equation}
\begin{split}
\sum_\sigma \M_\sigma(X) \SS^\sigma(A)
&=\sum_\sigma \M_\sigma(X)\sum_{\tau\le\sigma} \G_\tau(A) \\
&=\sum_\tau \left( \sum_{\sigma\ge\tau} \M_\sigma(X) \right)\G_\tau(A)\,. \\
\end{split}
\end{equation}
Equating the coefficients of $\G_\tau$ in both expressions, we see
that for any permutation $\tau$ with descent composition $I$, we must have
\begin{equation}
\label{sigtau}
\sum_{\sigma\ge\tau}\M_\sigma(X) = F_I(X) \,.
\end{equation}

\noindent
Now, given a permutation $\sigma$, the set of permutations greater than, or
equal to $\sigma$ in the left permutohedron order that are $\omega(I)$ for
some $I$ depends only on $C(\sigma)$ and is equal to
\begin{equation}
\{\omega(I) | I\succeq C(\sigma)\}.
\end{equation}
Indeed, it is the case if $\sigma$ is itself an $\omega(J)$ and one then only
has to see that all permutations $\tau$ greater than $\sigma$ satisfy
$C(\tau)\succeq C(\sigma)$.

Thanks to this property, we can assume by induction that the lemma is true for
all permutations strictly greater than $\sigma$ in the left permutohedron.
Then, if $\sigma = \omega(J)$, by inversion on the lattice of
compositions, $\M_{\omega(J)}(X)=M_J(X)$, as expected. Otherwise, $\sigma$ is
strictly greater than $\omega(C(\sigma))$, so that it has to be zero to
satisfy~(\ref{sigtau}).
\qed

Note that the multiplication formula \mref{eq:multS} is valid as soon as
the $\sigma_i$'s are such that
$\sigma=\sigma_1\carr\sigma_2\carr\cdots\carr\sigma_r$
(it not necessarily has to be the maximal factorisation). 
Also, from \mref{eq:multS}, one obtains the coproduct of a $\M_\sigma$,
as given in \cite{AS}, Theorem~3.1.

Dually, the multiplication formula for $\M_\alpha\M_\beta$ (Theorem 4.1
of~\cite{AS}) is equivalent to the computation of $\Delta \SS^\gamma$. Again,
it is more easily obtained on this side. To state the coproduct
formula, it will be convenient to introduce the notation
\begin{equation}
\check\SS^\sigma=\SS^{\sigma^{-1}}\,.
\end{equation}
Then,
\begin{equation}
\Delta\check\SS^\sigma=
\sum_{u\cdot v\le\sigma}\<\sigma|u\shuffle v\>
\check\SS^{\std(u)}\otimes\check\SS^{\std(v)}\,,
\end{equation}
that is, we sum over pairs of complementary subwords $u$, $v$ of $\sigma$
whose concatenation $u\cdot v$ is a permutation smaller than $\sigma$ in the
\emph{right} weak order.

For example,

\begin{equation}
\begin{split}
\Delta \SS^{32451}
&= 1 \otimes \SS^{32451} + \SS^{1} \otimes \SS^{2134}
   + 2 \SS^{1} \otimes \SS^{2341} + \SS^{21} \otimes \SS^{231} \\
 & + 2 \SS^{21} \otimes \SS^{123} + \SS^{321} \otimes \SS^{12}
   +   \SS^{213} \otimes \SS^{21} + \SS^{2134} \otimes \SS^{1} \\
 & +   \SS^{3241} \otimes \SS^{1} + \SS^{32451} \otimes    1.
\end{split}
\end{equation}

Indeed, here are the pairs of complementary subwords of $52134$ whose
concatenation $u\cdot v$ is smaller than $52134$:
\begin{equation}
\begin{split}
& (\epsilon,52134),\ (5,2134),\ (2,5134),\ (1,5234),\ (52,134),\ (51,234),\\
& (21,534),\ (521,34),\ (213,54),\ (5213,4),\ (2134,5),\ (52134,\epsilon).
\end{split}
\end{equation}

By the argument of \cite[Prop. 3.6]{NCSF6}, which holds for any multiplicative
basis, we have:
\begin{proposition}
The family $(\M_\sigma)$ where $\sigma$ goes along the set of
anticonnected permutations is a linear basis of the primitive Lie algebra of
$\FQSym$.
\end{proposition}

\subsection{Analog of the pair $(e,f)$ (first kind)} 

The basis $\Lambda^I$ of product of elementary symmetric functions of $\Sym$
can be extended to $\FQSym$ in the same way as the $S^I$. One has here to set
\begin{equation}
\LL^\sigma=\sum_{\tau\ge\sigma}\G_\tau\,,
\end{equation}
so that, for example
\begin{equation}
\LL^{132}=\G_{132} + \G_{231} + \G_{321}
\end{equation}
and
\begin{equation}
\LL^{3412}=\G_{3412} + \G_{3421} + \G_{4312} + \G_{4321}.
\end{equation}

\begin{proposition}
The basis $(\LL^\sigma)$ is multiplicative, meaning that
\begin{equation}
\LL^\sigma \LL^\tau = \LL^{\sigma \bullet \tau}\,.
\end{equation}  
\end{proposition}

\Proof The product rules of the basis $(\G_\sigma)$ can be written as
\begin{equation}
\G^\mu \G^\nu =
   \sum_{\alpha\in\SG_{m,n}} \G_{\alpha \circ (\mu\bullet\nu)}\,,
\end{equation}
where $\SG_{m,n}$ is set of the inverse of the permutations occurring in the
shuffle $(12\dots m)\shuffle(m+1\dots m+n)$. It follows that
\begin{equation}
\LL^\sigma \LL^\tau = 
   \sum_{\alpha\in\SG_{m,n}\,,\ \mu\ge\sigma\,,\ \nu\ge\tau}
   \G_{\alpha \circ (\mu\bullet\nu)}\,.
\end{equation}

As before, it is well-known that the decomposition
$\alpha \circ (\mu\bullet\nu) \in \SG_{m,n} \circ (\SG^m\times\SG^n)$
is reduced and that
\begin{equation}
\alpha \circ (\mu\bullet\nu) \ge \sigma \bullet \tau
\qquad\text{if and only if}\qquad
\mu\bullet\nu \ge \sigma \bullet \tau\,,
\end{equation}
which completes the proof. 
\qed

Then the noncommutative symmetric function $\Lambda^i= R_{(1^i)}$ is equal to
the free quasi-symmetric function $\LL^{i\,i-1\,\dots\,1}$.
Hence the correct identification is
\begin{equation}
\Lambda^I=\Lambda^{(i_1i_2\dots i_k)}=
\LL^{      (i_1\,i_1-1\,\dots\,1)\ 
  \bullet\ (i_1+i_2\,\dots\,i_1+1)\ \bullet \dots
  \bullet\ (n\,\dots\,i_1+i_2+\dots i_{k-1}+1) }
\end{equation}
that is
\begin{equation}
\Lambda^I=\LL^{\alpha(I^c)}=\LL^{\alpha(\bar I^\sim)}
\end{equation}
where $\alpha(I^c)$ is the shortest permutation whose descent set is the
complementary set of the descent set of $I$.

\bigskip
The adjoint basis $(\WW_{\mu})$ of $(\LL_{\sigma})$ is an analog of the
forgotten basis. It is given by
\begin{equation}
\F_\sigma = \sum_{\tau \le \sigma} \WW_\tau\,,
\end{equation}
or equivalently
\begin{equation}
\WW_\tau = \sum_{\sigma \le \tau} \mu(\tau,\sigma)\F_\sigma\,,
\end{equation}
where $\mu(\tau,\sigma)$ is the Moebius function of the
permutoedron~\cite{Bj,Ede}.
The same argument as in the previous section shows that

\begin{proposition}
(i) The family $(\WW_\sigma)$ where $\sigma$ goes along the set of connected
permutation is a linear basis of the primitive Lie algebra of $\FQSym$.
  
\noindent
(ii) The commutative image $\WW_\sigma(X)$ vanishes unless $\sigma$ is of the
form $\alpha(I)$ for a composition $I$.
\end{proposition}

\subsection{Analog of the pair $(e,f)$ (second kind)}

\subsubsection{Definition}

The Poirier-Reutenauer basis $\G^\sigma$, defined in~\cite{PR} by
\begin{equation}
\G^\sigma=\G_{\sigma_1}\G_{\sigma_2}\cdots \G_{\sigma_r}
\end{equation}
where $\sigma=\sigma_1\bullet\sigma_2\bullet\cdots\bullet\sigma_r$
is the factorization of $\sigma$ into connected permutations, has to be
considered as another analog of the basis of products of elementary symmetric
functions. Indeed, it contains the noncommutative elementary functions
$\Lambda_n=\G^{\omega_n}=\G_{\omega_n}$, where $\omega_n=n\cdots 1$ and, when
$\sigma=\alpha(I)$ is the minimum element of a descent class $D_I$, one has
\begin{equation}
\G^{\alpha(I)}=\Lambda^{I^c}\,.
\end{equation}
For example, with $I=(3,2,1)$, $\alpha(I)=124365$ and
\begin{equation}
\G^{124365}= \G_{1}\G_{1}\G_{21}\G_{21} = \Lambda^{1122}.
\end{equation}

It is convenient to introduce also the basis $\F^{\sigma}$, defined by
\begin{equation}
\F^\sigma=\F_{\sigma_1}\F_{\sigma_2}\cdots \F_{\sigma_r}
\end{equation}
Since the factorization into connected components commutes with inversion, one
has
\begin{equation}
\F^\sigma=\G^{\sigma^{-1}}\,.
\end{equation}

\subsubsection{Transition matrices and order relations}

To determine the transition matrices between the pairs of bases
$(\F^\sigma)$, $(\F_\sigma)$ and $(\G^\sigma)$, $(\G_\sigma)$, we introduce
two relations between permutations.
Given two permutations $\sigma$ and $\tau$, we write $\sigma\prec_F \tau$ iff
$\F_\tau$ appears in the expansion of $\F^\sigma$ in this basis.
We define in the same way $\prec_G$:
\begin{equation}\label{ordfg}
\sigma\prec_G \tau\Longleftrightarrow \sigma^{-1}\prec_F \tau^{-1}\,,
\end{equation}

By definition of the scalar product in $\FQSym$, the relations can
equivalently be defined by
\begin{equation}\label{ordf}
\sigma\prec_F \tau\Longleftrightarrow \< \F^\sigma\,,\,\G_\tau\>=1
\end{equation}
\begin{equation}\label{ordg}
\sigma\prec_G \tau\Longleftrightarrow \< \G^\sigma\,,\,\F_\tau\>=1.
\end{equation} 

\begin{lemma} 
The relations $\prec_F$ and $\prec_G$ are order relations.
\end{lemma}

\Proof
We can restrict ourselves to the $\prec_F$ case thanks to~(\ref{ordfg}).
By definition, $\prec_F$ is antisymmetric since the product
$\F_{\sigma_1}\F_{\sigma_2}$ consists in $\sigma_1\bullet\sigma_2$ plus other
terms greater in the lexicographic order.
In other words, the transition matrix between the bases $\F^\sigma$ and
$\F_\sigma$ is triangular with ones on the diagonal and this is equivalent to
the reflexivity and antisymmetry of $\prec_F$.

Now, let us assume $\sigma\prec_F\tau\prec_F\nu$. This means that
$\tau$ appears in the shuffle of the connected components factors of $\sigma$.
The same words for $\nu$ and $\tau$.
Since the connected components of $\tau$ are themselves shuffles of connected
components of $\sigma$, we deduce $\sigma\prec_F\nu$, hence proving the
transitivity of $\prec_F$.
\qed

\begin{figure}[ht]
\[\epsfig{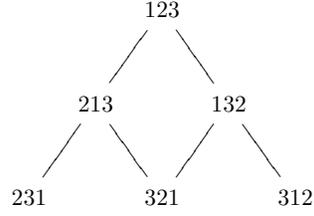}\]
\caption{\label{cat3}The catenohedron order for $n=3$.}
\end{figure}

\begin{figure}[ht]
\[\epsfig{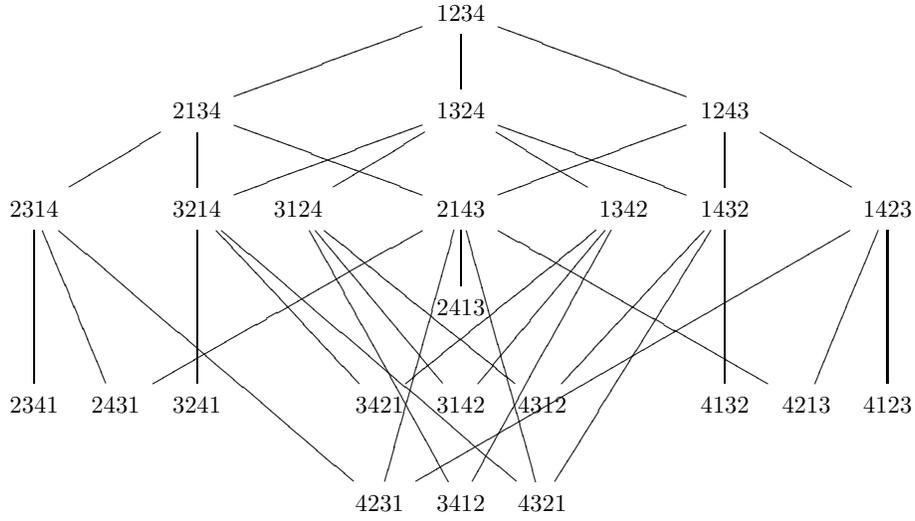}\]
\caption{\label{cat4}The catenohedron order for $n=4$.}
\end{figure}

The order $\prec_F$ will also be called the \emph{catenohedron order} (see
Figures~\ref{cat3} and~\ref{cat4}).
%

\begin{note}
The map $\sigma\ra \deg(\sigma):=n-l(J(\sigma))$ is strictly increasing.
\end{note}

If one denotes by $\mu_F$ et $\mu_G$ the Moebius functions of the above
orders, the transition matrices are
\begin{equation}
\label{passage}
\F^\sigma=\sum_{\sigma\prec_F \tau} \F_\tau;\quad
\F_\sigma=\sum_{\sigma\prec_F \tau}\mu_F(\sigma,\tau) \F^\tau\
\end{equation}
\begin{equation}
\G^\sigma=\sum_{\sigma\prec_G \tau} \G_\tau;\quad
\G_\sigma=\sum_{\sigma\prec_G \tau}\mu_G(\sigma,\tau) \G^\tau\,.
\end{equation}

Note that $\mu_G(\sigma,\tau)=\mu_F(\sigma^{-1},\tau^{-1})$. 

\begin{equation}
M_3(\F^\sigma,\F_\tau) =
\left(
\begin{matrix}
1 & . & . & . & . & . \\
1 & 1 & . & . & . & . \\
1 & . & 1 & . & . & . \\
1 & . & 1 & 1 & . & . \\
1 & 1 & . & . & 1 & . \\
1 & 1 & 1 & . & . & 1 
\end{matrix}
\right)
\end{equation}

\begin{equation}
M_3(\F_\tau,\F^\sigma) =
\left(
\begin{matrix}
1  & .  & .  & . & . & . \\
-1 & 1  & .  & . & . & . \\
-1 & .  & 1  & . & . & . \\
.  & .  & -1 & 1 & . & . \\
.  & -1 & .  & . & 1 & . \\
1  & -1 & -1 & . & . & 1 
\end{matrix}
\right)
\end{equation}

We shall now see that those matrices can be described in combinatorial terms. 

\subsubsection{Explicitation of $\prec_F$ and computation of $\mu_F$}
     
A chain of length $n$ from $\sigma$ to $\tau$ is a sequence of distinct
elements
\begin{equation}
\sigma=\sigma_0\prec_F\cdots \sigma_{n-1}\prec_F\sigma_n=\tau\,.
\end{equation}  
The set of such chains will be denoted by $\cat_n(\sigma,\tau)$. 
One has
\begin{equation}
\mu_F(\sigma,\tau) =
\sum_{n\rm\  odd} \# \cat_n(\sigma,\tau) -
\sum_{n\rm\ even} \# \cat_n(\sigma,\tau)
\end{equation} 
which proves that the restriction of $\mu_F$ to an interval of $\prec_F$
coincides with the Moebius function of this interval.
Let us therefore study the intervals $[\sigma,\nu]_F$ for $\sigma\prec_F\nu$.                                                   

\begin{lemma}
Let $\sigma\prec_F \nu$. Any two elements $\sigma_1,\sigma_2$ of
$[\sigma,\nu]_F$ admit a supremum $\sigma_1\vee \sigma_2$.
\end{lemma} 

\Proof
We already know that $\tau\mapsto J(\tau)$ is strictly decreasing.
Let $\sigma_1,\sigma_2\in [\sigma,\nu]_F$, and let $J_i=J(\sigma_i)$ be the
corresponding compositions.
For $J=(j_1,\ldots,j_r)$, let $\Des(J)$ be the associated subset of $[1,n-1]$.
We now introduce projectors $\tilde p_T$, which are essentially the
ones defined in~\cite[Eq.(46)]{NCSF6}, with indexation by subsets, and
without standardisation.
Let $T=\{t_1<t_2<\cdots t_k\}\subset [1,n-1]$, $\tau_j$ being the subword
of $\tau$ obtrained by restricting the alphabet to $[t_{j-1}+1,t_j]$ (with
$t_0=0,t_{k+1}=n$).
We set
\begin{equation}
\tp_T(\tau):=\tau_1\cdot \tau_2\cdots \tau_k\cdot \tau_{k+1}\,.
\end{equation}
For example,
\begin{center}
\begin{tabular}{c||c|c|c|c|c|c|c|c}
$T$           & $\{\}$ & $\{1\}$ & $\{2\}$ & $\{3\}$ & $\{1,2\}$ &
$\{1,3\}$ & $\{2,3\}$ & $\{1,2,3\}$\\
\hline
$\tp_T(4231)$  & $4231$ & $1423$  & $2143$  & $2314$  & $1243$    & 
$1234$    & $2134$    & $1234$\\
\end{tabular}
\end{center}
Since the alphabets of the $\tau_j$ are disjoint intervals, one has
\begin{equation}
\Des (J(\tp_T(\sigma)))\supseteq T
\end{equation}
and the set of those $T$ for which this is an equality is stable
under intersection, whence the result.
\qed

In the above example, $\Des (\tp_{\{1,3\}}(4231))=\{123\}\supset \{1,3\}$, 
and the subsets for which equality is realized are
\begin{equation}
\{\} \  \{1\} \  \{2\} \  \{3\} \  \{1,2\} \
\{1,3\} \  \{2,3\} \  \{1,2,3\}\,.
\end{equation}
This family is indeed closed under intersection (but not under union).
We shall need the notion of \emph{immediate sucessors of $\sigma$ in the
direction~$\nu$}.

\begin{definition}
For $\sigma\prec_F\tau$, we shall denote by $M_{\sigma\tau}$ the set of
minimal elements of $]\sigma,\tau]=[\sigma,\tau]-\{\sigma\}$.
\end{definition}
 
\begin{proposition}
Let $\tau\in [\sigma,\nu]$. The Moebius function $\mu_F(\sigma,\tau)$ is
given by
\begin{equation}
\mu_F(\sigma,\tau)= \left\{
\begin{array}{cl}
(-1)^{\# M_{\sigma\tau}} & \ \hbox{if} \
\tau = \sup(M_{\sigma\tau}), \vtr{2} \\
0  & \ \hbox{otherwise}\,.
\end{array}
\right.
\end{equation}
\end{proposition}

\Proof
Let $\chi,\delta$ be the characteristic and Dirac functions. One has
\begin{equation}
\label{eq1}
\chi([\sigma,\tau])=\delta_\sigma +
\chi(\bigcup_{\sigma_1\in M_{\sigma\tau}}[\sigma_1,\tau])
\end{equation}
since $]\sigma,\tau]=\cup_{\sigma_1\in M_{\sigma\tau}}[\sigma_1,\tau]$. 
By inclusion-exclusion,   
\begin{equation}
\label{eq2}
\chi(\bigcup_{\sigma_1\in M_{\sigma\tau}}[\sigma_1,\tau])=
\sum_{k=1}^{\# M_{\sigma\tau}} (-1)^{k-1}
\sum_{F\subset M_{\sigma\tau}\ |F|=k}
\chi(\bigcap_{\sigma_1\in F}[\sigma_1,\tau])\,.
\end{equation}
For $F\subset M_{\sigma\tau}$, let $\sigma_F=sup(F)$. Then, \mref{eq1}
and \mref{eq2} yield
\begin{equation}
\delta_\sigma=\sum_{F\subset M_{\sigma\tau}}
  (-1)^{|F|}\chi([\sigma_F,\tau])\,.
\end{equation}
The result comes from the fact that $\sigma_F$ appears only once in the sum,
since the immediate successors are associated with subsets differing by
exactly one element from the subset associated with $\sigma$.

\begin{lemma}
Let $\sigma\prec_F \nu$. Then, if $\sigma\not= \nu$, there exists $\tau\in
[\sigma,\nu]$ such that $l(J(\sigma))=l(J(\tau))+1$.
\end{lemma}
 
Now, by definition of the Moebius function of $[\sigma,\tau]$, one has
\begin{equation}
\delta_\sigma=\sum_{\sigma\prec\sigma_1\prec\tau}
\mu_F(\sigma,\sigma_1) \chi([\sigma_1,\tau])
\end{equation}
whence the result.
\qed

\subsubsection{Transition between $\V_\sigma$ and $\F_\sigma$}
 
Let us recall that the basis $\V_\sigma$ of $\FQSym$ defined in~\cite{NCSF6}
is the adjoint basis of $\G^\sigma$.
Since $\G_\sigma$ is the adjoint basis of $\F_\sigma$, it appears that
\begin{equation}
\label{G2V}
\V_{\tau}=\sum_{\sigma\prec_G \tau}\mu_G(\sigma,\tau)\F_\sigma \,.
\end{equation}
  
\subsection{Analog of the pair $(h,m)$ (second kind)}

There is finally another basis that can be defined in the same way as the
first pair $(h,m)$.

\subsubsection{Definition}

Let $\sigma$ be a permutation and let us consider its decomposition 
$\sigma=\sigma_1\carr\sigma_2\carr\cdots\carr\sigma_r$ into anticonnected
components.
Then
\begin{equation}
\H^\sigma := \G_{\sigma_1} \G_{\sigma_2}\dots \G_{\sigma_r}.
\end{equation}
For example,
\begin{equation}
\begin{split}
\H^{3421} =& \G_{12} \G_{1}\G_{1} \\
=&\ \G_{1234} + \G_{1243} + \G_{1324} + \G_{1342} + \G_{1423} + \G_{1432}\\
+&\ \G_{2314} + \G_{2341} + \G_{2413} + \G_{2431} + \G_{3412} + \G_{3421}.
\end{split}
\end{equation}

It is easy to check that
\begin{equation}
\H^\sigma =\overline{\G^{\bar\sigma}}
\end{equation}
where the line over $\G$ has the following meaning: expand $\G^{\bar\sigma}$
in the $\G_\tau$ basis and then take the mirror image of all indices.
For example, one checks:
\begin{equation}
\H^{3421} = \overline{\G^{1243}} = \overline{\G_{1}\G_{1}\G_{21}}.
\end{equation}

In particular, this proves that the transition matrix between $\H_\sigma$ and
$\G_\tau$ is equal to the transition matrix between $\G^\sigma$ and $\G_\tau$
up to the same reordering of the indices of both bases.

\subsubsection{Duality}

Let $\UU_\sigma$ be the dual basis of $\H^\sigma$.
The same computations as before get immediately the transition matrices
between the $\UU$ and the $\F$.

\section{The splitting formula}

The \emph{internal product} $*$ is defined on $\FQSym$ by
\begin{equation}
\G_\sigma * \G_\tau = 
\begin{cases} 
\G_{\tau\sigma} \ \mbox{\rm if $|\tau|=|\sigma|$} \\
0 \ \mbox{\rm otherwise}
\end{cases}
\end{equation}
where $\tau\sigma$ is the usual product of the symmetric group.

The fundamental property for computing with the internal product in $\Sym$ is
the following \emph{splitting formula}:
\begin{proposition}[\cite{GKLLRT}]
\label{mackey}
Let $F_1,F_2,\ldots,F_r,G \in \Sym$. Then, 
\begin{equation}
(F_1F_2\cdots F_r)*G = \mu_r \left[ (F_1\otimes\cdots\otimes F_r)* \Delta^r G
\right]
\end{equation}
where in the right-hand side, $\mu_r$ denotes the $r$-fold ordinary
multiplication and $*$ stands for the operation induced on $\Sym^{\otimes n}$
by $*$.
\end{proposition}                       

This formula can be extended to $\FQSym$ with the following assumptions:
$F_1,\ldots,F_r$ are arbitrary, but $G$ must belong to the Patras-Reutenauer
algebra ${\cal A}$, defined in~\cite{PR}.

Indeed, by~(\ref{prodGsI}), for $I=(i_1,\dots,i_r)$ and $\sigma_k\in\SG_{i_k}$,
we have
\begin{equation}
\label{GSI}
(\G_{\sigma_1} \dots \G_{\sigma_r}) * H =
(\G_{\sigma_1} \bullet \dots \bullet \G_{\sigma_r}) * S^I * H
\end{equation}

If $H\in\NCSF$, the classical splitting formula~(\ref{mackey}) applies, and
\begin{equation}
\label{SIG}
\begin{split}
S^I * H
&= \mu_r \left[ (S_{i_1}\otimes\cdots\otimes S_{i_r})*_r \Delta^r H\right] \\
&= \sum_{(H)} (S_{i_1}*H_{(1)}) \dots (S_{i_r}*H_{(r)}),
\end{split}
\end{equation}
using Sweedler's notation, so that the right-hand side of~(\ref{GSI}) can be
rewritten as
\begin{equation}
\label{SIG2}
\begin{split}
& (\G_{\sigma_1} \bullet \dots \bullet \G_{\sigma_r}) * 
\sum_{(H); \deg H_{(k)}=i_k} (S_{i_1}*H_{(1)}) \dots (S_{i_r}*H_{(r)}) \\
&=
(\G_{\sigma_1} \bullet \dots \bullet \G_{\sigma_r}) * 
\left(\sum_{(H); \deg H_{(k)}=i_k} (H_{(1)}) \bullet\dots\bullet
(H_{(r)})\right)
* S^I \\
&= 
\sum_{(H)} (\G_{\sigma_1}*H_{(1)}) \dots (\G_{\sigma_r}*H_{(r)}) \\
&= \mu_r \left[ (\G_{\sigma_1}\otimes\G_{\sigma_r}) *_r \Delta^r H \right].
\end{split}
\end{equation}

By linearity, this proves that~(\ref{mackey}) is valid for arbitrary $F_1$,
$\dots$, $F_r\in\FQSym$ and $H\in\NCSF$.

Now, as observed by Schocker~\cite[Thm. 2.3]{Scho}, an identity from
Garsia-Reutenauer~\cite[Thm. 2.1]{GR}, implies that~(\ref{SIG}) is valid for
any $H=\pi_1\dots\pi_s$, if the $\pi_k$ are in the image of
$Lie_{j_k}\in\FQSym$ with the embedding $\sigma\mapsto\G_\sigma$.
Such products $\pi_1\dots\pi_s$ span the Patras-Reutenauer algebra
$\A = \oplus \A_n$~\cite{PaR}.
This is a Hopf subalgebra of $\FQSym$, and each $\A_n$ is stable under the
internal product. Summarizing, we have

\begin{proposition}
The splitting formula~(\ref{mackey}) holds for arbitrary $F_i\in\FQSym$ and
$H$ in the Patras-Reutenauer algebra $\A$.
\end{proposition}

The property was known to Schocker~\cite{Scho2}. A different proof has
been recently published by Blessenohl~\cite{Ble}.

\end{document}